\documentclass[12pt]{amsart}
\usepackage{amsmath}
\usepackage{amssymb}
\usepackage{hyperref}
\usepackage{epsfig}
\usepackage{graphicx}
\numberwithin{equation}{section}
\usepackage{color}

\input xy
\xyoption{all}

\calclayout
\allowdisplaybreaks[3]

\theoremstyle{plain}
\newtheorem{prop}{Proposition}

\newtheorem{theo}[prop]{Theorem}

\theoremstyle{definition}

\newtheorem{rema}[prop]{Remark}

\def\ra{\rightarrow}

\def\cE{{\mathcal E}}
\def\cF{{\mathcal F}}

\def\cL{{\mathcal L}}

\def\cO{{\mathcal O}}

\def\mq{{\mathfrak q}}
\def\mr{{\mathfrak r}}

\def\bfE{{\mathbf E}}

\def\bP{{\mathbb P}}

\def\bZ{{\mathbb Z}}

\def\bF{{\mathbb F}}

\def\CH{\mathrm{CH}}

\def\Bl{\mathrm{Bl}}
\def\CH{\mathrm{CH}}

\def\lim{\mathrm{lim}}

\def\Sym{\mathrm{Sym}}

\makeatother
\makeatletter

\author{Brendan Hassett}
\address{Department of Mathematics\\
Brown University \\
Box 1917 
151 Thayer Street
Providence, RI 02912 \\
USA}
\email{bhassett@math.brown.edu}

\author{Alena Pirutka}
\address{Courant Institute\\
                New York University \\
                New York, NY 10012 \\
                USA }
\email{pirutka@cims.nyu.edu}

\author{Yuri Tschinkel}
\address{Courant Institute\\
                New York University \\
                New York, NY 10012 \\
                USA }
\email{tschinkel@cims.nyu.edu}

\address{Simons Foundation\\
160 Fifth Avenue\\
New York, NY 10010\\
USA}

\title[Quartic double fourfolds]{A very general quartic double fourfold is not stably rational}

\begin{document}
\date{\today}

\begin{abstract}
We prove that a very general double cover of the projective four-space, ramified in a quartic threefold, is not stably rational. 
\end{abstract}

\maketitle

\section{Introduction}

In this note we consider quartic double fourfolds, 
i.e., hypersurfaces $X_f$ in the weighted projective space $\bP(2,1,1,1,1,1)$, with homogeneous 
coordinates $(s,x,y,z,t,u)$, given by
a degree four equation of the form
\begin{equation}
\label{eqn:eqn}
s^2+ f(x,y,z,t,u) =0. 
\end{equation}
The failure of stable rationality for cyclic covers of projective
spaces has been considered by Voisin \cite{Voisin}, Beauville \cite{beau-6},
Colliot-Th\'el\`ene--Pirutka
\cite{ct-pir-cyclic}, and Okada \cite{okada}.
We work over an uncountable ground field $k$ of characteristic zero. 
Our main result is

\begin{theo}
\label{theo:main}
Let $f\in k[x,y,z,t,u]$ be a very general degree four form. Then $X_f$ is not stably rational. 
\end{theo}

This note is inspired by \cite{Beau-AG}, which used the new 
technique of the decomposition of the diagonal
\cite{Voisin,ct-pirutka,totaro-JAMS}. The main difficulty is to  
construct a special $X$ in the family \eqref{eqn:eqn} with following properties:
\begin{itemize}
\item[({\bf O})] Obstruction: the second unramified cohomology group $H^2_{nr}(X)$ (or another birational invariant) 
does not vanish,
\item[({\bf R})] Resolution:  
there exists a resolution of singularities $\beta:\tilde{X} \ra X$, such that the morphism $\beta$ 
is universally $\CH_0$-trivial,
\end{itemize}
(see, e.g., Sections 2 and 4 of \cite{HPT} for definitions). 
The verification of both properties for potential examples of $X$ is notoriously difficult. 
The example considered in \cite{Beau-AG} satisfies 
the second property, but not the first. 

Our main goal here is to produce an $X$ satisfying both.
We have a candidate example:
\begin{equation}
\label{eqn:our}
V: s^2+xyt^2+xzu^2+yz(x^2+y^2+z^2-2(xy+xz+yz))=0.
\end{equation}
The singular locus of $V$ is a connected curve, 
consisting of 4 components: two nodal cubics, a conic, and a line. 
How do we find this example? We may transform equation (\ref{eqn:our})
to
\begin{equation}
\label{eqn:old}
yzs^2+xzt_1^2+xyu_1^2+(x^2+y^2+z^2-2(xy+xz+yz))v_1^2=0.
\end{equation}
Precisely, we homogenize via an additional variable $v$,
multiply through by $yz$, and absorb the squares into the variables
$t_1,u_1,$ and $v_1$. The resulting equation gives a bidegree $(2,2)$ hypersurface
$$V' \subset \bP^2 \times \bP^3,$$
birational to $V$ via the coordinate changes.
In \cite{HPT} we proved that this $V'$ satisfies both properties (O) and (R). 
In particular, $V$ also satisfies (O), since unramified cohomology is a birational invariant. 

A direct verification of property (R) for this $V$ is possible (so that we could take $X=V$), but we found it more transparent  
to take an alternative approach, applying the specialization argument twice:
First we can specialize a very general $X_f$ to a quartic double fourfold $X$
which is singular along a line $\ell$ (contained in the ramification locus); we choose $X$ to be very general subject to this condition.
Then we show that the blowup morphism
$$
\beta:\tilde{X}:=\Bl_{\ell}(X)\ra X
$$ 
is universally $\CH_0$-trivial 
and that $\tilde{X}$ is smooth, i.e., $X$ satisfies (R).  
Furthermore, there exists a quadric bundle structure $\pi:\tilde{X}\ra \bP^2$, 
with degeneracy divisor a smooth octic curve.
In Section~\ref{sect:geometry} we analyze this geometry. 
We consider a degeneration of these quadric bundles to a fourfold $X'$ 
which is birational to $V'$, and thus satisfies (O). 
The singularities of $X'$  are similar to those considered in \cite{HPT};
the verification of the required property (R) for $X'$ is easier in this presentation. 
This is the content of 
Section~\ref{sect:sing}.
In Section~\ref{sect:wrapup} we give the argument for failure of stable
rationality of very general (\ref{eqn:eqn}).

\

\noindent {\bf Acknowledgments:} The first author was 
partially supported through NSF grant 1551514.

\section{Geometry of quartic double fourfolds}
\label{sect:geometry}

Let $X\ra \bP^4$ be a double fourfold, ramified along a quartic threefold $Y\subset \bP^4$. 
From the equation \eqref{eqn:eqn} we see that the quartic double fourfold $X$ is singular precisely 
along the singular locus of the quartic threefold $Y\subset \bP^4$ given by $f=0$. 

We will consider quartic threefolds $Y$ double along $\ell$. These form a linear
series of dimension 
$$\binom{8}{4}-5-12=53$$
and taking into account changes of coordinates---automorphisms of 
$\bP^4$ stabilizing $\ell$---we have $34$ free parameters.

Let $\beta:\tilde{X}\ra X$ be the blowup of $X$ along $\ell$. 
We will analyze its properties by embedding
it into natural bundles over $\bP^2$.

We start by blowing up $\ell$ in $\bP^4$.
Projection from $\ell$
gives a projective bundle structure
$$\varpi:\Bl_{\ell}(\bP^4) \ra \bP^2$$
where we may identify
$$\Bl_{\ell}(\bP^4)\simeq \bP(\cE), \quad 
\cE=\cO_{\bP^2}^{\oplus 2} \oplus \cO_{\bP^2}(-1).$$
Write $h$ for the hyperplane class on $\bP^2$ and its pullbacks
and $\xi$ for the first Chern class of $\cO_{\bP(\cE)}(1)$. 
Taking global sections
$$\cO_{\bP^2}^{\oplus 5} \twoheadrightarrow \cE^{\vee}$$
induces morphisms
$$\bP(\cE) \hookrightarrow \bP(\cO_{\bP^2}^{\oplus})\simeq \bP^4 \times \bP^2;$$
projecting onto the first factor gives the blow up.
Its exceptional divisor 
$$E\simeq \bP(\cO_{\bP^2}^{\oplus 2})\simeq \bP^1 \times \bP^2$$
has class $\xi-h$. 

Let $\tilde{Y} \subset \bP(\cE)$ denote the proper transform of $Y$, which has class
$$4\xi-2E=2\xi+2h.$$
Conversely, divisors in this linear series map to quartic hypersurfaces
in $\bP^4$ singular along $\ell$.
Since $2\xi+2h$ is very ample in $\bP(\cE)$ the generic such divisor is smooth.
The morphism $\varpi$ realizes $\tilde{Y}$ as a conic bundle over $\bP^2$;
its defining equation $q$ may also be interpreted as a section of the
vector bundle $\Sym^2(\cE^{\vee})(2h)$.
Let $\gamma:\tilde{Y} \ra Y$ denote the resulting resolution; its exceptional
divisor $F=\tilde{Y} \cap E$ is a divisor of bidegree $(0,2)$ in $E\simeq \bP^1 \times \bP^2$.
Hence $F\ra \ell$ is a trivial conic bundle and 
$\gamma$ is universally $\CH_0$-trivial.

Let $\tilde{X} \ra \bP(\cE)$ denote the double cover branched over $\tilde{Y}$, i.e.,
$s^2=q.$
This naturally sits in the projectization of an extension
$$0 \ra \cL \ra \cF \ra \cE \ra 0,$$
where $\cL$ is a line bundle. 
Note the natural maps
$$\Sym^2(\cE^{\vee}) \hookrightarrow \Sym^2(\cF^{\vee})
\twoheadrightarrow \cL^{-2},$$
and their twists
$$\Sym^2(\cE^{\vee})(2h) \hookrightarrow \Sym^2(\cF^{\vee})(2h)
\twoheadrightarrow \cL^{-2}(2h);$$
the last sheaf corresponds to the coordinate $s$. 
%{\color{red}I did not understand if $r$ was defined/were we use it? Fixed: used to be $r$ now $s$ --Y.}
Since we are over $\bP^2$
the extension above must split; furthermore, the coordinate
$s$ induces a trivialization
$$\cL^{-2}(2h)\simeq \cO_{\bP^2}.$$
Thus we conclude
$$\cF \simeq \cO_{\bP^2}(1) \oplus \cE \simeq 
\cO_{\bP^2}(1) \oplus \cO_{\bP^2}^{\oplus 2} \oplus \cO_{\bP^2}(-1).$$

The divisor $\tilde{X} \subset \bP(\cF)$ is generically smooth; 
let $\beta:\tilde{X} \ra X$ denote the induced resolution of $X$.
Its exceptional divisor is a double cover of $E$ branched over $F$,
isomorphic to a product of a smooth quadric surface with $\bP^1$.
(A double cover of $\bP^2$ branched along a conic curve is a smooth
quadric surface.) It follows that $\beta$ is universally $\CH_0$-trivial.

We summarize the key elements we will need:
\begin{prop}
\label{prop:ram}
Let $X\ra \bP^4$ be a double fourfold, ramified along a quartic threefold $Y\subset \bP^4$. 
Assume that $Y$ is singular along a line $\ell$ and
generic subject to this condition.   
Let $\beta:\tilde{X}\ra X$ be the blowup of $X$ along $\ell$. 
Then $\tilde{X}$ is smooth and $\beta$ universally $\CH_0$-trivial. 
\end{prop}

Regarding $\tilde{X} \subset \bP(\cF)$, there is an induced quadric surface
fibration
$$\pi:\tilde{X} \ra \bP^2.$$
Let $D$ denote the degeneracy curve, naturally a divisor in
$$\det(\cF^{\vee}(2h))\simeq \cO_{\bP^2}(8).$$
The analysis above gives an explicit determinantal description of
the defining equation of $D$. Choose homogeneous forms
$$c \in \Gamma(\cO_{\bP^2}), \
F_1,F_2,F_3 \in \Gamma(\cO_{\bP^2}(2)), \
G_1,G_2 \in \Gamma(\cO_{\bP^2}(3)),
H \in \Gamma(\cO_{\bP^2}(4))$$
so that the symmetric matrix associated with $\tilde{X}$ takes the form:
$$\left( \begin{matrix} c  & 0 & 0 & 0 \\
                        0  & F_1 & F_2 & G_1 \\
                        0  & F_2 & F_3 & G_2 \\
                        0 &  G_1 & G_2 & H  
\end{matrix}\right)
$$

We fix coordinates to obtain a concrete equation for $\tilde{X}$.
Let $(x,y,z)$ denote coordinates of $\bP^2$, or equivalently,
linear forms on $\bP^4$ vanishing along $\ell$.
Let $s$ denote a local coordinate trivializing $\cO_{\bP^1}(1)\subset \cF$,
$t$ and $u$ coordinates corresponding to $\cO_{\bP^1}^{\oplus 2} \subset \cF$,
and $v$ to $\cO_{\bP^1}(-1)\subset \cF$. Then we have
\begin{equation} \label{eqn:general}
\tilde{X}=\{c s^2 + F_1t^2+ 2F_2tu + F_3u^2 + 2G_1tv + 2G_2uv + Hv^2=0\},
\end{equation}
where $F_1,F_2,F_3,G_1,G_2,$ and $H$ are homogeneous in $x,y,z$.

Finally, we interpret the degeneration curve in geometric
terms. Ignoring the constant, we may write
$$D={(F_1F_3-F_2^2)H-F_3G_1^2 + 2F_2G_1G_2-F_1G_2^2=0}.$$
Modulo $F_1F_3-F_2^2$ we have
$$-F_3G_1^2 + 2F_2G_1G_2-F_1G_2^2=0$$
which is equal to
$$\frac{-1}{F_1}(F_2G_1-F_1G_2)^2=\frac{-1}{F_3}(F_3G_1-F_2G_2)^2.$$

Thus we conclude that $D$ is tangent to a quartic plane curve
$$C=\{F_1F_3-F_2^2\}=0$$
at $16$ points. {\em Every} smooth quartic plane curve admits multiple
such representations: Surfaces
$$\{a^2F_1+2abF_2+b^2F_3=0\} \subset \bP^1_{a,b} \times \bP^2$$
are precisely degree two del Pezzo surfaces equipped with a 
conic bundle structure, the conic structures indexed by non-trivial
two-torsion points of the branch curve $C$.
One last parameter check: The moduli space of pairs
$(C,D)$ consisting of a plane quartic and a plane octic tangent at
$16$ points depends on
$$14+44-16-8=34$$
parameters. This is compatible with our first parameter count. 

\begin{rema}
{\em Smooth} divisors $\tilde{X} \subset \bP(\cF)$ as above 
necessarily have trivial Brauer group. This follows from
Pirutka's analysis \cite{pirutka-survol}: if the degeneracy
curve is smooth and irreducible then there cannot be 
unramified second cohomology. It also follows from a
singular version of the Lefschetz hyperplane theorem.
Let $\zeta=c_1(\cO_{\bP(\cF)}(1))$ so that $[\tilde{X}]=2\zeta+2h$.
This is almost ample: the line bundle $\zeta+h$
contracts the distinguished section $s:\bP^2 \ra \bP(\cF)$ associated with the
summand $\cO_{\bP^1}(1) \subset \cF$ to a point but otherwise induces
an isomorphism onto its image. In particular, $\zeta+h$ induces
a small contraction in the sense of intersection homology. 
The homology version of
the Lefschetz Theorem of Goresky-MacPherson \cite[p.~150]{GM} 
implies that $0\simeq H^3(\bP(\cF),\bZ) \stackrel{\sim}{\ra} H^3({\tilde X},\bZ)$.
\end{rema}

\section{Singularities of the special fiber}
\label{sect:sing}

We specialize (\ref{eqn:general}) to:
\begin{equation}
\label{eqn:special}
s^2+xyt^2+xzu^2+yz(x^2+y^2+z^2-2(xy+xz+yz))v^2=0.
\end{equation}

\begin{prop}
\label{prop:resolve}
The fourfold $X'\subset \bP(\mathcal F)$ defined by (\ref{eqn:special})
admits a resolution of singularities $\beta': \tilde{X}'\ra X'$ such that $\beta'$
is universally $\CH_0$-trivial. 
\end{prop}
The remainder of this section is a proof of this result.

\subsection{The singular locus}
\label{sect:sing-loc}

A direct computation in {\tt Magma} (or an analysis as in \cite[Section 5]{HPT})
yields that the singular locus of \eqref{eqn:special} 
is a connected curve consisting of the following components:
\begin{itemize}
\item 
Singular cubics:
\begin{align*}
E_z:=&\{ v^2 y(y-x)^2+ u^2 x = z = s = t = 0 \} \\
E_y:=&\{ v^2 z(z-x)^2+ t^2  x = y = s = u = 0 \} 
\end{align*}
\item 
Conics:
\begin{align*}
R_x:=&\{u^2 - 4v^2 + t^2 = x = z - y = s = 0 \} \\
C_x:=&\{zu^2 + yt^2 = s = v = x = 0 \} 
\end{align*}
\end{itemize}
The nodes of $E_z$ and $E_y$ are 
\begin{align*}
\mathfrak{n}_z:=&\{z=s=t=y-x=u=0\} \\
\mathfrak{n}_y:=&\{y=s=u=z-x=t=0\},
\end{align*}
respectively. Here $R_x$ and $C_x$ intersect transversally at two points,
$$
\mr_{\pm}:=\{u\pm it=v=s=z-y=x=0\};
$$ 
$R_x$ is disjoint from $E_z$ and $E_y$, 
and the other curves intersect transversally in a single point (in coordinates $(x,y,z)\times (s,t,u,v)$):
\begin{align*}
E_z\cap E_y=&\mq_x:=(1,0,0)\times (0,0,0,1),\\
E_z\cap C_x= &\mq_y:=(0,1,0)\times (0,0,1,0),\\
E_y\cap C_x=&\mq_z:=(0,0,1)\times (0,1,0,0).
\end{align*}
%Note that $E_z$, $E_y$, and  $C_x$ project, under $\pi:X'\ra \bP^2$ 
%to the coordinate lines $L_z,L_y$, and $L_x$, respectively. 
This configuration of curves is similar to the one considered in \cite{HPT}, but the singularities are different.

\subsection{Local \'etale description of the singularities and
resolutions}
\label{subsect:loc}
The structural properties of the resolution become clearer after
identifying \'etale normal forms for the singularities.

The main normal form is 
\begin{equation} \label{eqn:normal1}
a^2+b^2+c^2=p^2q^2
\end{equation}
which is singular along the locus
$$\{a=b=c=p=0\} \cup \{a=b=c=q=0\}.$$
This is resolved by successively blowing up along these 
components in either order. Indeed, after blowing up the 
first component, using $\{A,B,C,P\}$ for homogeneous 
coordinates associated with the corresponding generators
of the ideal, we obtain
$$A^2+B^2+C^2=P^2q^2.$$
The exceptional fibers are isomorphic to a non-singular
quadric hypersurface (when $q\neq 0$) or a quadric cone
(over $q=0$). Dehomogenizing by setting $P=1$, we obtain
$$A^2+B^2+C^2=q^2$$
which is resolved by blowing up $\{A=B=C=q=0\}$.
This has ordinary threefold double points at each point,
so the exceptional fibers are all isomorphic to non-singular
quadric hypersurfaces.

There are cases where 
$$\{a=b=c=p=0\} \cup \{a=b=c=q=0\}$$
are two branches of the same curve. For example,
this could arise from
\begin{equation}
\label{eqn:normal2}
a^2+b^2+c^2=(m^2-n^2-n^3)^2
\end{equation}
by setting $p=m-n\sqrt{1+n}$ and $q=m+n\sqrt{1+n}$.
Of course, we cannot pick one branch to blow up first.
We therefore blow up the origin first, using
homogeneous coordinates $A,B,C,D,P,Q$ corresponding
to the generators to obtain
$$A^2+B^2+C^2=P^2q^2=Q^2p^2.$$
The resulting fourfold is singular along the stratum
$$A=B=C=q=p=0$$
as well as the proper transforms of the original branches.
Indeed, on dehomogenizing $P=1$ we obtain local affine equation
$$A^2+B^2+C^2=Q^2p^2;$$
this is singular along $\{A=B=C=p=0\}$, the locus where the
exceptional divisor is singular, and $\{A=B=C=Q=0\}$, and proper transform
of $\{a=b=c=q=0\}$. 
The local affine equation is the same as (\ref{eqn:normal1}); we resolve by 
blowing up the singular locus of the exceptional divisor
followed by blowing up the proper transforms of the branches. This descends
to a resolution of (\ref{eqn:normal2}).

\subsection{Computation in local charts}
We exploit the symmetry under the involution exchanging $y\leftrightarrow z$
and $t \leftrightarrow u$. It suffices then to analyze $E_z,C_x,$ and $R_x$
and the distinguished points $\mathfrak{n}_z$, $\mq_x$, $\mq_y$, and $\mr_{+}$.

\subsubsection*{Analysis along the curve $C_x$}
Recall the equation of $X'$:
$$s^2+xyt^2+xzu^2+yz(x^2+y^2+z^2-2xy-2xz-2yz)v^2=0$$
and the equation of $C_x$:
$zu^2+yt^2=s=v=x=0$.
We order coordinates $(x,y,z), (s,t,u,v)$ and write intersections
\begin{itemize}
\item  $C_x\cap R_x=(0,1,1) \times (0,1,\pm i, 0)$;
\item  $C_x\cap E_z=(0,1,0) \times (0,0,1,0)$;
\item  $C_x\cap E_y=(0,0,1) \times (0,1,0,0)$.
\end{itemize}
We use the symmetry between $t$ and $u$ to reduce the number of cases.
 
\

\paragraph{{\it Chart} $u=1$, $z=1$}
We extract equations for the exceptional divisor $\bfE$ obtained
by blowing up $C_x$. In this chart,
$C_x$ takes the form
$$
1+yt^2=s=v=x=0
$$ and $X'$ is
$$s^2+x(yt^2+1)+v^2(y-1)^2+v^2xG=0,$$
where $v^2xG$ are the `higher order terms'.

Now we analyse the local charts of the blow up:
\begin{enumerate}
\item $\bfE: yt^2+1=0$,  $s=s_1(yt^2+1), x=x_1(yt^2+1), v=v_1(yt^2+1)$, the equation for $X'$, up to removing higher order terms, in new coordinates is:
$$s_1^2+x_1+v_1^2(y-1)^2=0,$$ this is smooth and rational.
The exceptional divisor 
$$s_1^2+x_1+v_1^2(y-1)^2=0, yt^2+1=0$$ is rational.
\item $\bfE: x=0$, $s=s_1x, v=v_1x, yt^2+1=wx$, equation of $X'$:
$$s_1^2+w+v_1(y-1)^2=0, yt^2+1=wx,$$ smooth;
\item $\bfE: s=0$, $x=x_1s, v=v_1s, yt^2+1=ws$:
$$1+wx_1+v_1(y-1)^2=0, yt^2+1=sw,$$ smooth.
\item $\bfE: v=0, s=s_1v, x=x_1v, yt^2+1=wv$, equation of $X'$ is
$$s_1^2+wx_1+(y-1)^2=0, yt^2+1=wv,$$ 
which has at most ordinary double singularity (corresponding to $C_x\cap R_x=\mr_{\pm}$) of type
$$a^2+b^2+cd=0, \  a=b=c=d=0.$$ 
This is resolved by one blowup.
\end{enumerate}

\paragraph{Chart $u=1$, $y=1$}
In this chart $C_x$ is
$z+t^2=s=v=x=0$ and $X'$ is
$$s^2+x(t^2+z)+v^2(z-1)^2+v^2xG=0,$$
where $v^2xG$ are the 'higher order terms'.
We analyze local charts of the blow up:
\begin{enumerate}
\item $\bfE: t^2+z=0$,  $s=s_1(t^2+z), x=x_1(t^2+z), v=v_1(t^2+z)$, the equation for $X'$, up to removing higher order terms, in new coordinates is:
$$s_1^2+x_1+v_1^2(z-1)^2=0,$$ this is smooth and rational.
The exceptional divisor 
$$s_1^2+x_1+v_1^2(z-1)^2=0, t^2+z=0$$ is rational.
\item $\bfE: x=0$, $s=s_1x, v=v_1x, t^2+z=wx$, equation of $X'$:
$$s_1^2+w+v_1(z-1)^2=0, t^2+z=wx,$$ smooth;
\item $\bfE: s=0$, $x=x_1s, v=v_1s, t^2+z=ws$:
$$1+wx_1+v_1(z-1)^2=0, t^2+z=sw,$$ smooth.
\item $\bfE: v=0, s=s_1v, x=x_1v, t^2+z=vw$, equation of $X'$ is
$$s_1^2+wx_1+(z-1)^2=0, t^2+z=wv,$$ 
or, up to removing the higher order terms
$$s_1^2+wx_1+(t^2+1)^2=0, z=-t^2+wv,$$
this has at most ordinary double singularities  
$$
s_1=w=x_1=0, t=\pm i
$$
(where we meet the proper transform of $R_x$) of type
$$a^2+b^2+cd=0, a=b=c=d$$
resolved as above by one blowup.
\end{enumerate}

\subsubsection*{Analysis near $\mathfrak{n}_z$}
Center the coordinates by setting $\xi=y-1$ 
$$s^2+(\xi+1)t^2+zu^2+(\xi+1)z(\xi^2+z^2-2z(\xi+2))=0.$$
Note that $E_z$ is given by
$$(\xi+1)\xi^2+u^2=z=s=t=0.$$
We regroup terms
$$s^2+(\xi+1)t^2+z(u^2+(\xi+1)\xi^2)+(\xi+1)z^2(z-2\xi-4)=0.$$
Provided $\xi\neq -1,-2$ this is \'etale-locally equal to
$$s_1^2+t_1^2+z_1(u^2+(\xi+1)\xi^2)+z_1^2=0$$
which is equivalent to normal form (\ref{eqn:normal2}). 
When $\xi=-1$ we are at the point $\mq_x$, which we analyze below.
A local computation at $\xi=-2$ shows that the singularity is
resolved there by blowing up $E_z$ and the exceptional fiber there
is isomorphic to $\bF_0$. In other words, we have ordinary 
threefold double points there as well.  

\subsubsection*{Blowing up the singular point $\mathfrak n_z$ of $E_z$}
  
The point $\mathfrak n_z$  lies in the chart $x=1, v=1$, where we now make computations.
The equation of the point (and the locus we blow up) is
$$s=t=u=z=y-1=0.$$
The equation of $X$  can be written as:
$$s^2+yt^2-2z^2y(y+1)+zu^2+z^3y+yz(y-1)^2=0.$$
The curve $E_z$ has equations
$$y(y-1)^2+u^2=z=s=t=0.$$
Now we compute the charts for the blow up:
  
\begin{enumerate}
\item 
$\bfE: s=0$. The change of variables is  $u=su_1, t=st_1, z=sz_1, y=1+y_1s.$ Then the equation of $X'$ (resp. the exceptional divisor $\bfE$), up to removing the higher order terms, is:
  $$1+t_1^2(1+y_1s)-2z_1^2(1+sy_1)(2+sy_1)=0$$
  (resp. $1+t_1^2-2z_1^2=0$), so that the blow up and the exceptional divisor  are smooth, and $\bfE$ is rational.
\item 
$\bfE: t=0$. The change of  variables is  $s=s_1t,  u=u_1t, z=z_1t, y=1+y_1t:$ the equations are
  $$s_1^2+(1+y_1t)-2z_1^2(1+y_1t)(2+y_1t)=0,$$  and $E$ is given by
  $$s_1^2+1-4z_1^2=0,$$ so that the blowup is  smooth at any point of the exceptional divisor.
\item 
$\bfE: z=0$, the change of variables is $s=s_1z, u=u_1z,  y=1+y_1z$; we obtain
$$
s_1^2+(1+y_1z)t_1^2-2(1+y_1z)(2+y_1z)=0
$$ 
and the equation of $\bfE$ is $s_1^2+t_1^2-4=0, $ so that the blow up is smooth at any point of the exceptional divisor.
\item 
$\bfE: y_1:=y-1=0$, the change of variables is $z=z_1y_1, s=s_1y_1, u=u_1y_1, t=t_1y_1$; the equations are
  $$s_1^2+t_1^2(1+y_1)-2z_1^2(1+y_1)(2+y_1)+u_1^2y_1z+z_1^3y_1(1+y_1)+z_1(1+y_1)=0,$$ 
this is smooth, as well as the exceptional divisor ($y_1=0$).
  
\item $\bfE: u=0$, the change of variables is $s=s_1u, t=t_1u, z=z_1u, y=1+y_1u $, 
the equations for the proper transform of $X'$ are
$$
s_1^2+(1+y_1u)t_1^2-2z_1^2(1+y_1u)(2+y_1u)+uz_1+uz_1^3(1+y_1u)+z_1uy^2_1(1+y_1u)=0,
$$
and the proper transform of $E_z$ is given by
$$
(1+y_1u)y_1^2+1=z_1=s_1=t_1=0.
$$
The exceptional divisor 
$$\bfE: s_1^2+t_1^2-4z_1^2=0$$
is singular along $s_1=t_1=z_1=u=0$ (and $y_1$ is free). The
resulting curve is denoted $R_z\simeq \bP^1$; note that $R_z$ meets
the proper transform of $E_z$ at two points $y_1=\pm i$.
  \end{enumerate}
  
\subsubsection*{Blowing up $R_z$}
For the analysis of singularities we can remove higher order terms, so that the equation of the variety (resp. $R_z$) is given by:
  $$s_1^2+t_1^2-4z_1^2+uz_1+uz_1y^2_1=0$$
  and $s_1=t_1=z_1=u=0$.
  
The charts for the new blow up with exceptional divisor $\bfE'$ are:
  
\begin{enumerate}
\item $\bfE': s_1=0$, then after the usual change of variables for a blow up, we obtain 
the equation:
$$
1+t_2^2-4z_2^2+u_2z_2+u_2z_2y^2_1=0,
$$ 
which is smooth.
  \item $\bfE': t_1=0$ is similar to the previous case.
  \item $\bfE': z_1=0$, we obtain the equation 
  $$s_2^2+t_2^2-4+u_2+u_2y_1^2=0,$$ that is smooth;
  \item $\bfE': u=0$,  we obtain the equation 
  $$s_2^2+t_2^2-4z_2^2+z_2(1+y^2_1)=0,$$ 
which has ordinary double points at $s_2=t_2=z_2=y_1^2+1=0$.
These are resolved by blowing up the proper transform of $E_z$.
  \end{enumerate}

\subsubsection*{Analysis near $\mq_x$}
Dehomogenize
$$
s^2+xyt^2+xzu^2+yz(x^2+y^2+z^2-2(xy+xz+yz))v^2=0
$$ 
by setting $v=1$ and $x=1$ to obtain
$$
s^2+yt^2+zu^2+yz(1+y^2+z^2-2(y+z+yz))=0.
$$
We first analyze at $\mq_x$, the origin in this coordinate system.
Note that $1+y^2+z^2-2(y+z+yz) \neq 0$ here and thus its
square root can be absorbed (\'etale locally) in to $s,t,$ and $u$
to obtain
$$s_1^2+yt_1^2+zu_1^2+yz=0.
$$
Setting $y_1=y+u_1^2$ and $z_1=z+t_1^2$ gives
$$s_1^2+y_1z_1=t_1^2u_1^2,$$
which is equivalent to normal form (\ref{eqn:normal1}).
(The blow up over the generic point of $E_z$ was analyzed previously.)

\subsubsection*{Blowing up $R_x$}

Similar to the analysis of singularities near $R_z$, see also \cite[Section 5.2 (4)]{HPT}

\subsection{Summary of the resolution}
\label{sect:summary}
\subsubsection*{Blowup steps}
The resolution $\beta'$ is a sequence of blowups:
\begin{enumerate}
\item{Blow up the nodes $\mathfrak{n}_z$ and $\mathfrak{n}_y$; the resulting
fourfold is singular along rational curves $R_z$ and $R_y$ in the exceptional locus,
meeting the proper transforms of $E_z$ and $E_y$ transversally in two points
sitting over $\mathfrak{n}_z$ and $\mathfrak{n}_y$, respectively.}
\item{The exceptional divisors are quadric threefolds singular along $R_z$ and $R_y$.}
\item{At this stage the singular locus consists of six smooth rational
curves, the proper
transforms of $E_z,E_y,R_x,C_x$ and the new curves $R_z$ and $R_y$, with a total
of nine nodes. (This is the configuration appearing in \cite[Section 5]{HPT}.)}
\item{The local analytic structure is precisely as indicated in Section~\ref{subsect:loc}.
Thus we can blow up the six curves in any order to obtain a resolution of singularities.
The fibers are either the Hirzebruch surface $\bF_0$ or a union of
Hirzebruch surfaces $\bF_0\cup_{\Sigma}\bF_2$ where $\Sigma\simeq \bP^1$
with self intersections $\Sigma^2_{\bF_0}=2$ and $\Sigma^2_{\bF_2}=-2$.}
\end{enumerate}
For concreteness, we blowup in the order
$$R_z,R_y,E_z,E_y,C_x,R_x.$$

\subsubsection*{Exceptional fibers}
The following fibers arise:
\begin{itemize}
\item{Over the nodes $\mathfrak{n}_z$ and $\mathfrak{n}_y$:
The exceptional fiber has two three-dimensional
components. One is the standard resolution of a 
quadric threefold singular along a line, that is,
$$
{\mathbf F}'=\bP(\cO_{\bP^1}^{\oplus 2} \oplus \cO_{\bP^1}(-2)).
$$
The other is a quadric surface fibration $\mathbf{F}''\ra \bP^1$,
over $R_z$ and $R_y$ respectively,
smooth except for two fibers corresponding to
the intersections with $E_z$ and $E_y$; the singular fibers are unions 
$\bF_0 \cup \bF_2$ as indicated above. The intersection 
$\mathbf{F}'\cap \mathbf{F}''$ is along the distinguished subbundle
$$\bP(\cO_{\bP^1}^{\oplus 2}) \subset \mathbf{F}'$$
which meets the smooth fibers of $\mathbf{F}''\ra \bP^1$ in 
hyperplanes and the singular fibers in smooth rational
curves in $\bF_2$ with self-intersection $2$.}
\item{Over $E_z$: the exceptional divisor is a quadric
surface fibration over $\bP^1$, with two degenerate fibers
of the form $\bF_0 \cup \bF_2$ corresponding to the intersections
with $C_x$ and $E_y$.}
\item{Over $E_y$: the exceptional divisor is a quadric
surface fibration with one degenerate fiber, corresponding to the
intersection with $C_x$.}
\item{Over $C_x$: a quadric surface fibration with two degenerate 
fibers corresponding to the intersections with $R_x$.}
\item{Over $R_x$: a smooth quadric surface fibration.}
\end{itemize}
In each case, the fibers of $\beta'$ are universally
$\CH_0$-trivial.

\section{Proof of the Theorem~\ref{theo:main}}
\label{sect:wrapup}

We recall implications of 
the ``integral decomposition 
of the diagonal and specialization'' method, following 
\cite{ct-pirutka},  \cite{Voisin}.

\begin{theo}
\label{thm:main-help}
\cite[Theorem 2.1]{Voisin}, \cite[Theorem 1.14 and Theorem 2.3]{ct-pirutka}
Let 
$$
\phi: \mathcal X\ra B
$$ 
be a flat projective morphism  of complex varieties with smooth generic fiber. 
Assume that there exists a point $b\in B$ so that the fiber 
$$
X:=\phi^{-1}(b)
$$ 
satisfies the following conditions:
\begin{itemize}
\item $X$ admits a desingularization 
$$
\beta: \tilde{X}\ra X,
$$
where the morphism $\beta$ is universally $\CH_0$-trivial,
\item $\tilde{X}$ is not universally $\CH_0$-trivial. 
\end{itemize}
Then a very general fiber of $\phi$ is not universally $\CH_0$-trivial and, in particular, not
stably rational.
\end{theo}

We apply this twice: 
Consider a family of double fourfolds $X_f$ ramified along a quartic threefold $f=0$, as in (\ref{eqn:eqn}). Let $X'$ be the fourfold given by (\ref{eqn:special}) and let $V'$ be the bidegree $(2,2)$ hypersurface defined in (\ref{eqn:old}).

\begin{itemize}
\item [(1)]  As mentioned in the introduction, $V'$ satisfies
property (O); this is an application of Pirutka's computation
of unramified second
cohomology of quadric surface bundles over $\bP^2$ \cite{pirutka-survol}. 
By construction, $X'$ is birational to $V'$. 
Proposition~\ref{prop:resolve} and Section~\ref{sect:summary} yield property (R) for $X'$. We conclude that
very general hypersurfaces $\tilde{X} \subset \bP(\cF)$ given by Equation \ref{eqn:general}, in
Section~\ref{sect:geometry}, following Proposition~\ref{prop:ram}, 
fail to be universally $\CH_0$-trivial.
\item [(2)]
 By
Proposition~\ref{prop:ram},  the resolution morphism  $\beta:\tilde{X} \ra X$ is universally $\CH_0$-trivial; here $X$ is a double fourfold, ramified along a quartic which is singular along a line.
A second application of Theorem~\ref{thm:main-help} 
to the family of double fourfolds ramified along a quartic threefold
completes the proof of Theorem~\ref{theo:main}. 
\end{itemize}

\bibliographystyle{alpha}
\bibliography{double}
\end{document}